\documentclass{amsproc}

\usepackage{amssymb}

\newtheorem{theorem}{Theorem}[section]
\newtheorem{lemma}[theorem]{Lemma}
\newtheorem{proposition}[theorem]{Proposition}

\theoremstyle{definition}

\newtheorem{example}[theorem]{Example}

\newtheorem{question}[theorem]{Question}
\newtheorem{notation}[theorem]{Notation}
\theoremstyle{remark}

\numberwithin{equation}{section}

\newcommand{\abs}[1]{\lvert#1\rvert}

\newcommand{\R}{{\mathbb R}}
\newcommand{\Z}{{\mathbb Z}}
\newcommand{\FF}{{\mathcal F}}
\newcommand{\GG}{{\mathcal G}}

\newcommand{\C}{{\mathbb C}}

\newcommand{\PP}{{\mathbb P}}

\newcommand{\HH}{{\mathbb H}}
\newcommand{\PPP}{{\mathcal P}}
\newcommand{\DD}{{\mathbb D}}

\begin{document}
\title[Equidistribution theorem for harmonic measures]
{Weak form of equidistribution theorem for harmonic measures of
foliations by hyperbolic surfaces}

{}
\author{Shigenori Matsumoto}
\address{Department of Mathematics, College of
Science and Technology, Nihon University, 1-8-14 Kanda, Surugadai,
Chiyoda-ku, Tokyo, 101-8308 Japan
}
\email{matsumo@math.cst.nihon-u.ac.jp
}
\thanks{The author is partially supported by Grant-in-Aid for
Scientific Research (C) No.\ 25400096.}


\date{\today }
\begin{abstract}
We show that the equidistribution theorem of C. Bonatti and
 X. G\'omez-Mont for a special kind of foliations by hyperbolic surfaces
does not hold in general, and seek for a weaker form
valid for general foliations  by hyperbolic surfaces.
\end{abstract}

\maketitle

\section{Introduction}

Let $M$ be a smooth closed manifold, and let $\FF$ be a smooth foliation
by hyperbolic surfaces, i.e.\ a 2 dimensional foliation equipped with
a smooth leafwise metric $g_P$ of constant curvature $-1$. 
Let 
$v_P$ be the leafwise Poincar\'e volume form, and for a point $z\in M$ and
$\rho>0$, let $B_\rho(z)$ be the leafwise $\rho$-disk centered at $z$.
When $B_\rho(z)$ is an embedded disk in $M$, let $\beta_\rho(z)$ be the probability measure of $M$ supported on
$B_\rho(z)$
defined by
$$
\beta_\rho(z)=\frac{1}{\int_{B_\rho(z)}v_P}\,v_P\vert_{B_\rho(z)}.$$
When $B_\rho(z)$ is not embedded, define 
$\beta_\rho(z)$ using the universal cover of the leaf.

In \cite{BGM}, Christian Bonatti and Xavier G\'omez-Mont has shown the following
theorem.

\begin{theorem}
 Let $\Sigma$ be a closed oriented hyperbolic surface, and let
$\Phi:\pi_1(\Sigma)\to PSL(2,\C)$ be a nonelementary representation.
Endow leaves of  the associated foliated $\PP^1$ bundle $(N,\GG)$ 
with a hyperbolic metric lifted from $\Sigma$.
Then there exists a probability measure $\mu$ on $N$ such that
for any sequences $z_n\in N$ and $\rho_n\to\infty$,
 $\beta_{\rho_n}(z_n)$ converges weakly to $\mu$.
\end{theorem}

The measure $\mu$ turns out to be the unique harmonic measure
of the foliation $\GG$ in the sense of \cite{G}. See Section 4 for
more detail.
 Thus one may ask the following question.

\begin{question}
 For $(M,\FF)$ as above, if $\beta_{\rho_n}(z_n)$ converges weakly to a
 measure $\mu$ as $\rho_n\to\infty$, is it true that $\mu$
is a harmonic measure of $\FF$?
\end{question}

In Section 2, we shall answer this question in the negative, and in
Section 3, we propose a measure $\mu_{\rho,\rho'}(z)$ modified for
the positive answer. In Section 4, we raise a further question and
give a new example
 of foliations for which the conclusion of the theorem
of Bonatti and G\'omez-Mont holds.

\section{A counterexample}

Let $Solv_3$ be the 3-dimensional unimodular solvable nonnilpotent Lie
group. The multiplication of $Solv_3=\{(x,q,t)\}$ is given by
$$
(x,q,t)(x',q',t')=(e^tx'+x,\,e^{-t}q'+q,\,t+t').$$
It has  a structure of semidirect product:
$$1\to\R^2\to Solv_3\to\R\to 1.$$
Any lattice $\Gamma$ of $Solv_3$ is a semidirect product
$$1\to\Z^2\to\Gamma\to\Z\to 1$$
such that $\Gamma\cap\R^2=\Z^2$.
The multiplication is given by
$$
(n,m,\ell)(n',m',\ell')=((n',m')A^\ell+(n,m),\,\ell+\ell')$$
for some hyperbolic matrix $A\in SL(2,\Z)$.
The quotient manifold $M=\Gamma\setminus Solv_3$ is a $T^2$ bundle over
$S^1$ with monodromy $^t\!\!A$:
$$
T^2=\Z^2\setminus\R^2\to M=\Gamma\setminus Solv_3\to S^1=\Z\setminus\R.$$
Denote by
$G=\{q=0\}$ the subgroup of $Solv_3$, isomorphic to the
2-dimensional solvable nonabelian Lie group, and let $\tilde\FF$ be the
orbit foliation of the right $G$ action. 
Notice that the leaf passing
through $(x_0,q_0,t_0)$ is just $L_{q_0}=\{q=q_0\}$. 
The left action of the lattice $\Gamma$ commutes with the right $G$
action,
and therefore $\tilde\FF$ descends to a foliation $\FF$
on $M$.
Now 
$$g=e^{-2t}dx^2+e^{2t}dq^2+dt^2$$
is a left invariant metric on $Solv_3$.
The restriction of  $g$ to each leaf $L_{q}$ of $\tilde \FF$
is written as
$$g_P=e^{-2t}dx^2+dt^2.$$
If we change the variable by $y=e^t$, then we get
$$
g_P=(dx^2+dy^2)/y^2,$$
the Poincar\'e metric on the half plane $\HH$. 
That is, we have an identification 
$$Solv_3\ni(x,q,t)\leftrightarrow (x+e^ti,q)\in\HH\times\R,$$
where the right $G$ action leaves each leaf $L_q=\HH\times\{q\}$ invariant.
The action of the one parameter subgroup 
$\{Y^t=(0,0,t)\}$ of $G$ on each leaf $L_{q}\cong\HH$ is given by
$$Y^t(x+yi)=x+e^tyi,$$
 and the one parameter subgroup $\{S^s=(s,0,0)\}$
by 
$$S^s(x+yi)=ys+x+yi.$$
They satisfy
$$
Y^t\circ S^s=S^{se^{-t}}\circ Y^t.$$
See Figure 1.
\begin{figure}[h]
{\unitlength 0.1in%
\begin{picture}( 27.3900, 13.7200)( 11.9000,-19.9200)%
%
\special{pn 8}%
\special{pa 1377 803}%
\special{pa 3738 803}%
\special{fp}%
%
\special{pn 8}%
\special{pa 1381 1384}%
\special{pa 3720 1384}%
\special{fp}%
\special{pa 1372 1667}%
\special{pa 3743 1667}%
\special{fp}%
%
\special{pn 20}%
\special{pa 1190 1964}%
\special{pa 3929 1964}%
\special{fp}%
%
\special{pn 8}%
\special{pa 1654 620}%
\special{pa 1654 1984}%
\special{fp}%
\special{pa 2564 620}%
\special{pa 2564 1992}%
\special{fp}%
\special{pa 3465 620}%
\special{pa 3465 1960}%
\special{fp}%
%
\special{pn 20}%
\special{pa 1654 1355}%
\special{pa 1654 845}%
\special{fp}%
\special{sh 1}%
\special{pa 1654 845}%
\special{pa 1634 912}%
\special{pa 1654 898}%
\special{pa 1674 912}%
\special{pa 1654 845}%
\special{fp}%
%
\special{pn 20}%
\special{pa 1654 1643}%
\special{pa 1654 1427}%
\special{fp}%
\special{sh 1}%
\special{pa 1654 1427}%
\special{pa 1634 1494}%
\special{pa 1654 1480}%
\special{pa 1674 1494}%
\special{pa 1654 1427}%
\special{fp}%
%
\special{pn 20}%
\special{pa 2564 1652}%
\special{pa 2564 1436}%
\special{fp}%
\special{sh 1}%
\special{pa 2564 1436}%
\special{pa 2544 1503}%
\special{pa 2564 1489}%
\special{pa 2584 1503}%
\special{pa 2564 1436}%
\special{fp}%
%
\special{pn 20}%
\special{pa 2573 1331}%
\special{pa 2573 856}%
\special{fp}%
\special{sh 1}%
\special{pa 2573 856}%
\special{pa 2553 923}%
\special{pa 2573 909}%
\special{pa 2593 923}%
\special{pa 2573 856}%
\special{fp}%
%
\special{pn 20}%
\special{pa 3465 1648}%
\special{pa 3465 1445}%
\special{fp}%
\special{sh 1}%
\special{pa 3465 1445}%
\special{pa 3445 1512}%
\special{pa 3465 1498}%
\special{pa 3485 1512}%
\special{pa 3465 1445}%
\special{fp}%
%
\special{pn 20}%
\special{pa 3465 1340}%
\special{pa 3465 864}%
\special{fp}%
\special{sh 1}%
\special{pa 3465 864}%
\special{pa 3445 931}%
\special{pa 3465 917}%
\special{pa 3485 931}%
\special{pa 3465 864}%
\special{fp}%
%
\special{pn 20}%
\special{pa 1731 797}%
\special{pa 2510 797}%
\special{fp}%
\special{sh 1}%
\special{pa 2510 797}%
\special{pa 2443 777}%
\special{pa 2457 797}%
\special{pa 2443 817}%
\special{pa 2510 797}%
\special{fp}%
%
\special{pn 20}%
\special{pa 2655 808}%
\special{pa 3401 808}%
\special{fp}%
\special{sh 1}%
\special{pa 3401 808}%
\special{pa 3334 788}%
\special{pa 3348 808}%
\special{pa 3334 828}%
\special{pa 3401 808}%
\special{fp}%
%
\special{pn 20}%
\special{pa 1741 1384}%
\special{pa 2160 1384}%
\special{fp}%
\special{sh 1}%
\special{pa 2160 1384}%
\special{pa 2093 1364}%
\special{pa 2107 1384}%
\special{pa 2093 1404}%
\special{pa 2160 1384}%
\special{fp}%
%
\special{pn 20}%
\special{pa 2646 1384}%
\special{pa 3056 1384}%
\special{fp}%
\special{sh 1}%
\special{pa 3056 1384}%
\special{pa 2989 1364}%
\special{pa 3003 1384}%
\special{pa 2989 1404}%
\special{pa 3056 1384}%
\special{fp}%
%
\special{pn 20}%
\special{pa 1749 1667}%
\special{pa 1936 1667}%
\special{fp}%
\special{sh 1}%
\special{pa 1936 1667}%
\special{pa 1869 1647}%
\special{pa 1883 1667}%
\special{pa 1869 1687}%
\special{pa 1936 1667}%
\special{fp}%
%
\special{pn 20}%
\special{pa 2651 1667}%
\special{pa 2864 1667}%
\special{fp}%
\special{sh 1}%
\special{pa 2864 1667}%
\special{pa 2797 1647}%
\special{pa 2811 1667}%
\special{pa 2797 1687}%
\special{pa 2864 1667}%
\special{fp}%
\end{picture}}%

\caption{$Y^t$ moves points upwards along the geodesics by length $t$.
$S^s$ moves points horizontally along the horocycles by length
 $s$. $Y^t$ contracts $S^s$.}
 \end{figure}
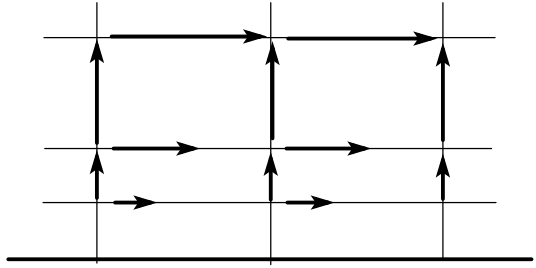
On the other hand, the left $Solv_3$ action (in particular $\Gamma$
action) leaves the product structure invariant, and the action on the first
factor $\HH$ is given by
$$(x,q,t)\cdot z=e^tz+x.$$
It is
not only $g_P$-isometric but also 
leaves the point $\infty$ on $\partial \HH$ invariant.
That is, each leaf of the foliation $\FF$ of
the quotient manifold $M$ admits a pointed hyperbolic structure.

The flow $\{S^s\}$ leaves the coordinate $y$, whence the old coordinate
$t$, invariant. Thus it leaves fibers of the fibration $T^2\to M\to S^1$
invariant and is a linear flow on it parallel to an eigenvector of
the matrix $^t\!A$.

Now $m=dx\wedge dq\wedge dt$ is a biinvariant Haar measure of $Solv_3$. 
If we denote by $v_P=y^{-2}dx\wedge dy$ the leafwise Poincar\'e volume form
of $\tilde\FF$, then 
$$
m=-y\,v_P\wedge dq.$$ 
The measure $m$ yields a probability measure on $M$, also denoted by $m$.
By a criterion in \cite{G},  $m$ is a harmonic measure of $\FF$, 
since the function $y$
is a harmonic function on $\HH$. Moreover by a general theorem
of Bertrand Deroin and Victor Kleptsyn \cite{DK}, it is the unique
harmonic measure.
The rest of this section is devoted to the proof of the following
theorem.

\begin{theorem}
 There exist $z_n\in M$ such that $\beta_n(z_n)$ converges to $\mu\neq m$.
\end{theorem}

We consider an infinite cyclic covering $\hat M$ of $M$ and the lift
$\hat\FF$ of the foliation $\FF$. Precisely,
$$
\hat M=\Z^2\setminus Solv_3=T^2\times\R,$$
where $\Z^2$ is the normal subgroup of the lattice $\Gamma$.
Let us denote by $\PPP(\hat M)$ the space of the Radon probability measures
of $\hat M$, endowed with the pointwise convergence topology
on the space $C_0(\hat M)$ of continuous functions on $\hat M$
with compact support.

Every leaf of $\hat\FF$ is pointedly isometric to $\HH$.
Choose one leaf and identify it with $\HH$.
For $\rho>0$, let $z_\rho=e^\rho i\in\HH\subset\hat M$. Notice that
the hyperbolic distance of $z_\rho$ to the horocycle $\{y=1\}$ is
$\rho$. We shall show that the probability measure
$\beta_\rho(z_\rho)\in\PPP(\hat M)$
converges to a measure $\hat\mu\in\PPP(\hat M)$ as $\rho\to\infty$.
The boundary $\partial B_\rho(z_\rho)$ of the disk $B_\rho(z_\rho)$
is tangent to $\{y=1\}$ and satisfies the equation:
$$
x^2+(y-\frac{R+1}{2}\,)^2=\frac{1}{4}(R-1)^2, \ \mbox{ where }\
R=e^{2\rho}.$$
See Figure 2.
\begin{figure}[h] 
{\unitlength 0.1in%
\begin{picture}( 26.8000, 16.1800)( 8.7000,-16.7000)%
%
\special{pn 8}%
\special{ar 2008 1003 535 535  0.0000000  6.2831853}%
%
\special{pn 8}%
\special{pa 1374 1540}%
\special{pa 2823 1540}%
\special{fp}%
%
\special{pn 20}%
\special{pa 1378 1725}%
\special{pa 2811 1725}%
\special{fp}%
%
\special{pn 8}%
\special{pa 2000 465}%
\special{pa 2819 465}%
\special{dt 0.045}%
%
\special{pn 8}%
\special{pa 1370 902}%
\special{pa 2844 902}%
\special{fp}%
%
\special{pn 8}%
\special{pa 1387 1171}%
\special{pa 2848 1171}%
\special{fp}%
%
\special{pn 4}%
\special{pa 2130 902}%
\special{pa 1861 1171}%
\special{fp}%
\special{pa 2105 902}%
\special{pa 1836 1171}%
\special{fp}%
\special{pa 2080 902}%
\special{pa 1811 1171}%
\special{fp}%
\special{pa 2055 902}%
\special{pa 1786 1171}%
\special{fp}%
\special{pa 2029 902}%
\special{pa 1761 1171}%
\special{fp}%
\special{pa 2004 902}%
\special{pa 1735 1171}%
\special{fp}%
\special{pa 1979 902}%
\special{pa 1710 1171}%
\special{fp}%
\special{pa 1954 902}%
\special{pa 1685 1171}%
\special{fp}%
\special{pa 1929 902}%
\special{pa 1660 1171}%
\special{fp}%
\special{pa 1903 902}%
\special{pa 1635 1171}%
\special{fp}%
\special{pa 1878 902}%
\special{pa 1609 1171}%
\special{fp}%
\special{pa 1853 902}%
\special{pa 1584 1171}%
\special{fp}%
\special{pa 1828 902}%
\special{pa 1559 1171}%
\special{fp}%
\special{pa 1803 902}%
\special{pa 1534 1171}%
\special{fp}%
\special{pa 1777 902}%
\special{pa 1509 1171}%
\special{fp}%
\special{pa 1752 902}%
\special{pa 1496 1158}%
\special{fp}%
\special{pa 1727 902}%
\special{pa 1492 1137}%
\special{fp}%
\special{pa 1702 902}%
\special{pa 1488 1116}%
\special{fp}%
\special{pa 1677 902}%
\special{pa 1483 1095}%
\special{fp}%
\special{pa 1651 902}%
\special{pa 1479 1074}%
\special{fp}%
\special{pa 1626 902}%
\special{pa 1479 1049}%
\special{fp}%
\special{pa 1601 902}%
\special{pa 1475 1028}%
\special{fp}%
\special{pa 1576 902}%
\special{pa 1475 1003}%
\special{fp}%
\special{pa 1551 902}%
\special{pa 1475 978}%
\special{fp}%
\special{pa 1525 902}%
\special{pa 1479 948}%
\special{fp}%
\special{pa 1500 902}%
\special{pa 1483 919}%
\special{fp}%
\special{pa 2155 902}%
\special{pa 1887 1171}%
\special{fp}%
\special{pa 2181 902}%
\special{pa 1912 1171}%
\special{fp}%
\special{pa 2206 902}%
\special{pa 1937 1171}%
\special{fp}%
\special{pa 2231 902}%
\special{pa 1962 1171}%
\special{fp}%
\special{pa 2256 902}%
\special{pa 1987 1171}%
\special{fp}%
\special{pa 2281 902}%
\special{pa 2013 1171}%
\special{fp}%
\special{pa 2307 902}%
\special{pa 2038 1171}%
\special{fp}%
\special{pa 2332 902}%
\special{pa 2063 1171}%
\special{fp}%
\special{pa 2357 902}%
\special{pa 2088 1171}%
\special{fp}%
\special{pa 2382 902}%
\special{pa 2113 1171}%
\special{fp}%
\special{pa 2407 902}%
\special{pa 2139 1171}%
\special{fp}%
\special{pa 2433 902}%
\special{pa 2164 1171}%
\special{fp}%
\special{pa 2458 902}%
\special{pa 2189 1171}%
\special{fp}%
\special{pa 2483 902}%
\special{pa 2214 1171}%
\special{fp}%
\special{pa 2508 902}%
\special{pa 2239 1171}%
\special{fp}%
\special{pa 2529 906}%
\special{pa 2265 1171}%
\special{fp}%
\special{pa 2533 927}%
\special{pa 2290 1171}%
\special{fp}%
\special{pa 2538 948}%
\special{pa 2315 1171}%
\special{fp}%
\special{pa 2542 969}%
\special{pa 2340 1171}%
\special{fp}%
\special{pa 2542 994}%
\special{pa 2365 1171}%
\special{fp}%
\special{pa 2542 1020}%
\special{pa 2391 1171}%
\special{fp}%
\special{pa 2538 1049}%
\special{pa 2416 1171}%
\special{fp}%
\special{pa 2538 1074}%
\special{pa 2441 1171}%
\special{fp}%
\special{pa 2529 1108}%
\special{pa 2466 1171}%
\special{fp}%
%
\special{pn 4}%
\special{pa 2521 1141}%
\special{pa 2491 1171}%
\special{fp}%
\put(28.8200,-9.3600){\makebox(0,0)[lb]{$e^{t+\Delta t}$}}%
\put(28.9500,-11.9600){\makebox(0,0)[lb]{$e^t$}}%
\put(28.8600,-15.7400){\makebox(0,0)[lb]{$1$}}%
\put(28.8600,-4.8200){\makebox(0,0)[lb]{$R=e^{2\rho}$}}%
%
\put(40.5000,-21.7000){\makebox(0,0)[lb]{}}%
%
\special{pn 4}%
\special{sh 1}%
\special{ar 1980 1300 16 16 0  6.28318530717959E+0000}%
\special{sh 1}%
\special{ar 1980 1300 16 16 0  6.28318530717959E+0000}%
\put(20.5000,-13.9000){\makebox(0,0)[lb]{$z_\rho$}}%
\end{picture}}%

\caption{}
\end{figure}
Putting $y=e^t$, we get
$$
x=\pm\sqrt{R(e^t-1)+e^t-e^{2t}}.$$
Now since $v_P=y^{-2}dx\wedge dy=e^{-t}dx\wedge dt$, the area
$A(R,t,\Delta t)$ of the set $B_\rho(z_\rho)\cap(T^2\times[t,t+\Delta t])$
is given by
$$
A(R,t,\Delta t)=
2\int_t^{t+\Delta t}\sqrt{R(e^{-t}-e^{-2t})+e^{-t}-1}\,\,dt.$$
On the other hand, the area $A(R)$ of $B_\rho(z_\rho)$ is given by
$$
A(R)=\pi(e^\rho+e^{-\rho}-2)=\pi(R^{1/2}+R^{-1/2}-2).$$
Now we have
$$
\beta_\rho(z_\rho)(T^2\times[t,t+\Delta t])
=\frac{A(R,t,\Delta t)}{A(R)}
=\frac{2}{\pi}\int_t^{t+\Delta
t}\frac{\sqrt{R(e^{-t}-e^{-2t})+e^{-t}-1}}
{R^{1/2}+R^{-1/2}-2}dt.$$
Therefore the limit measure $\hat\mu$ as $\rho\to\infty$ should satisfy
$$
\hat\mu(T^2\times[t,t+\Delta t])=\frac{2}{\pi}\int_t^{t+\Delta t}\sqrt{e^{-t}-e^{-2t}}dt.$$
On the other hand, the portion of the measure $\beta_\rho(z_\rho)$
supported on $T^2\times[t,t+\Delta t]$ ($t>0$) becomes more and more
invariant by the flow $S^s$ as $\rho\to\infty$, 
while $S^s$ is a linear flow of irrational slope on $T^2\times\{t\}$
and $dx\wedge dq$ is the unique measure invariant by $S^s$. Therefore one concludes
that
$$\beta_\rho(z_\rho)\to\hat\mu=\hat\Phi\, dx\wedge dq\wedge dt
\mbox{ as }\rho\to\infty,$$
where
$$\displaystyle
\hat\Phi(t)=\left\{\begin{array}{ll}
\frac{2}{\pi}\sqrt{e^{-t}-e^{-2t}} & \mbox{ if }t\geq 0, \\
                  0  & \mbox{ if }t\leq0.
		   \end{array}\right.$$
The actual proof needs the evaluation on a function from 
$C_0(\hat
M)$, which is a routine and omitted.
But $\beta_\rho(z_\rho)\to\hat\mu$ does not garantee that
$\hat\mu$ is a probability measure, since the constant 1
does not belong to $C_0(\hat M)$ and some part of $\beta_\rho(z_\rho)$
may escape to $\infty$. This, however, is assured by the following concrete
computation:
$$
\int_0^\infty\sqrt{e^{-t}-e^{-2t}}dt=\pi/2.$$
Also this implies a stronger fact that $\beta_\rho(z_\rho)\to\hat\mu$ pointwise
on any bounded continuous function.
The function $\hat\Phi$ takes the maximum value at $t=\log2$.
See Figure 3. 
\begin{figure}[h]
{\unitlength 0.1in%
\begin{picture}( 36.6400, 11.5200)( 18.2000,-16.2900)%
%
\special{pn 8}%
\special{pa 1820 1389}%
\special{pa 5476 1389}%
\special{fp}%
%
\special{pn 8}%
\special{ar 3548 1389 1328 830  3.1415927  3.8495597}%
%
\special{pn 8}%
\special{pa 2548 844}%
\special{pa 2580 834}%
\special{pa 2612 823}%
\special{pa 2643 813}%
\special{pa 2675 804}%
\special{pa 2707 796}%
\special{pa 2738 789}%
\special{pa 2769 784}%
\special{pa 2801 780}%
\special{pa 2832 778}%
\special{pa 2863 778}%
\special{pa 2893 780}%
\special{pa 2924 785}%
\special{pa 2954 792}%
\special{pa 2984 800}%
\special{pa 3014 811}%
\special{pa 3044 823}%
\special{pa 3104 849}%
\special{pa 3133 864}%
\special{pa 3163 879}%
\special{pa 3193 893}%
\special{pa 3223 908}%
\special{pa 3252 923}%
\special{pa 3312 951}%
\special{pa 3342 964}%
\special{pa 3372 978}%
\special{pa 3402 991}%
\special{pa 3432 1003}%
\special{pa 3462 1016}%
\special{pa 3582 1064}%
\special{pa 3613 1075}%
\special{pa 3673 1097}%
\special{pa 3704 1107}%
\special{pa 3734 1118}%
\special{pa 3764 1128}%
\special{pa 3795 1137}%
\special{pa 3825 1147}%
\special{pa 3887 1165}%
\special{pa 3917 1174}%
\special{pa 3948 1183}%
\special{pa 4010 1199}%
\special{pa 4040 1207}%
\special{pa 4071 1214}%
\special{pa 4102 1222}%
\special{pa 4164 1236}%
\special{pa 4195 1242}%
\special{pa 4226 1249}%
\special{pa 4258 1255}%
\special{pa 4320 1267}%
\special{pa 4351 1272}%
\special{pa 4383 1277}%
\special{pa 4414 1282}%
\special{pa 4446 1287}%
\special{pa 4477 1292}%
\special{pa 4509 1296}%
\special{pa 4540 1300}%
\special{pa 4572 1304}%
\special{pa 4603 1308}%
\special{pa 4635 1312}%
\special{pa 4731 1321}%
\special{pa 4762 1324}%
\special{pa 4794 1326}%
\special{pa 4826 1329}%
\special{pa 4986 1339}%
\special{pa 5050 1341}%
\special{pa 5082 1343}%
\special{pa 5115 1344}%
\special{pa 5179 1346}%
\special{pa 5211 1346}%
\special{pa 5275 1348}%
\special{pa 5339 1348}%
\special{pa 5371 1349}%
\special{pa 5484 1349}%
\special{fp}%
%
\special{pn 8}%
\special{pa 2860 784}%
\special{pa 2860 1389}%
\special{dt 0.045}%
\put(27.7000,-15.5000){\makebox(0,0)[lb]{$\log2$}}%
\put(21.7000,-15.2000){\makebox(0,0)[lb]{$0$}}%
\end{picture}}%

\caption{The function $\hat\Phi$.}
\end{figure}
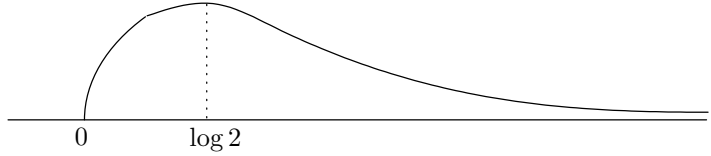
Returning to the compact manifold $M$, the previous observation
shows that the limit measure $\mu$ of
$\beta_\rho(z_\rho)$ is the projected image of $\hat\mu$ and is written
as $\mu=(\Phi\circ p)\, m$, where $p:M\to S^1$ is the bundle projection
and $\Phi$ is a continuous function on
$S^1$ given by 
$$\Phi(t)=\sum_{k\in\Z}\hat\Phi(t+k).$$
But $\Phi$ is not a constant function
since $\Phi(0)<\Phi(\log2)$, showing that $\mu\neq m$, as is required.

The author is greatful to Hiroki Kodama for showing him this simple proof.

\section{Weak form of equidistribution}

Let $\DD=\{z\in\C\mid\abs{z}<1\}$ be the disk model of the Poincar\'e plane. 
For $0<R<1$, denote  $\DD(R)=\{\abs{z}<R\}$ and let $\rho$ be
the Poincar\'e distance from $0$ to the circle $\partial\DD(R)$, i.e.\
$$
\displaystyle \rho=\frac{1}{2}\log\frac{1+R}{1-R}.$$
For $R<R'<1$, let
$\rho'$ be the Poincar\'e distance between $\partial\DD(R)$ and
$\partial\DD(R')$. Define a function
$\psi_{\rho,\rho'}:\DD\to[0,\infty)$ by
$$
\psi_{\rho,\rho'}(z)=\left\{ \begin{array}{ll} 1 & \mbox{\ if \ }\abs{z}\leq R
\\ \frac{R'-\abs{z}}{R'-R} & \mbox{\ if \ } R\leq\abs{z}\leq R'\\
0 & \mbox{\ if \ } R'\leq\abs{z}<1.\end{array}
\right.$$ 
The function $\psi_{\rho,\rho'}$ is
determined by $\rho$ and $\rho'$.
See Figure 4.
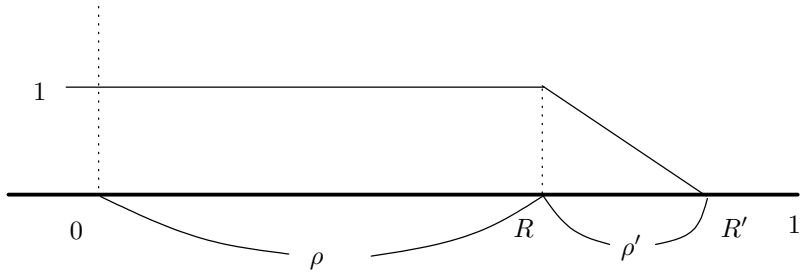
\begin{figure}[h]
{\unitlength 0.1in%
\begin{picture}( 41.3100, 13.0700)( 10.4000,-18.6700)%
%
\special{pn 20}%
\special{pa 1040 1545}%
\special{pa 5171 1545}%
\special{fp}%
%
\special{pn 8}%
\special{pa 1512 560}%
\special{pa 1512 1535}%
\special{dt 0.045}%
%
\special{pn 8}%
\special{pa 1342 982}%
\special{pa 3844 982}%
\special{fp}%
%
\special{pn 8}%
\special{pa 3834 992}%
\special{pa 3834 1545}%
\special{dt 0.045}%
%
\special{pn 8}%
\special{pa 1522 1555}%
\special{pa 1609 1597}%
\special{pa 1639 1610}%
\special{pa 1668 1624}%
\special{pa 1697 1637}%
\special{pa 1726 1651}%
\special{pa 1756 1664}%
\special{pa 1785 1677}%
\special{pa 1815 1689}%
\special{pa 1844 1701}%
\special{pa 1904 1725}%
\special{pa 1934 1736}%
\special{pa 2024 1766}%
\special{pa 2086 1784}%
\special{pa 2116 1791}%
\special{pa 2147 1799}%
\special{pa 2179 1806}%
\special{pa 2241 1818}%
\special{pa 2273 1823}%
\special{pa 2304 1828}%
\special{pa 2368 1838}%
\special{pa 2400 1842}%
\special{pa 2431 1846}%
\special{pa 2463 1850}%
\special{pa 2495 1855}%
\special{pa 2507 1856}%
\special{fp}%
%
\special{pn 8}%
\special{pa 2939 1867}%
\special{pa 3099 1842}%
\special{pa 3195 1824}%
\special{pa 3226 1818}%
\special{pa 3258 1811}%
\special{pa 3289 1804}%
\special{pa 3351 1788}%
\special{pa 3382 1779}%
\special{pa 3412 1770}%
\special{pa 3442 1760}%
\special{pa 3502 1738}%
\special{pa 3531 1726}%
\special{pa 3560 1713}%
\special{pa 3589 1699}%
\special{pa 3617 1685}%
\special{pa 3645 1670}%
\special{pa 3701 1638}%
\special{pa 3728 1621}%
\special{pa 3756 1605}%
\special{pa 3783 1587}%
\special{pa 3810 1570}%
\special{pa 3834 1555}%
\special{fp}%
%
\special{pn 8}%
\special{pa 3834 1545}%
\special{pa 3872 1599}%
\special{pa 3892 1626}%
\special{pa 3912 1651}%
\special{pa 3933 1675}%
\special{pa 3956 1697}%
\special{pa 3979 1718}%
\special{pa 4005 1736}%
\special{pa 4032 1751}%
\special{pa 4060 1765}%
\special{pa 4090 1777}%
\special{pa 4121 1787}%
\special{pa 4152 1796}%
\special{pa 4184 1805}%
\special{pa 4186 1806}%
\special{fp}%
%
\special{pn 8}%
\special{pa 4427 1806}%
\special{pa 4461 1799}%
\special{pa 4495 1791}%
\special{pa 4527 1782}%
\special{pa 4558 1771}%
\special{pa 4586 1758}%
\special{pa 4611 1741}%
\special{pa 4633 1721}%
\special{pa 4650 1696}%
\special{pa 4665 1669}%
\special{pa 4676 1638}%
\special{pa 4687 1606}%
\special{pa 4696 1573}%
\special{pa 4698 1565}%
\special{fp}%
\put(13.6200,-17.7600){\makebox(0,0)[lb]{$0$}}%
\put(36.8300,-17.6600){\makebox(0,0)[lb]{$R$}}%
\put(47.6900,-17.6600){\makebox(0,0)[lb]{$R'$}}%
\put(26.1800,-19.3700){\makebox(0,0)[lb]{$\rho$}}%
\put(42.4600,-18.8700){\makebox(0,0)[lb]{$\rho'$}}%
\put(51.2000,-17.4600){\makebox(0,0)[lb]{$1$}}%
%
\special{pn 8}%
\special{pa 3834 975}%
\special{pa 4692 1552}%
\special{fp}%
%
\put(35.2600,-13.1000){\makebox(0,0)[lb]{}}%
\put(11.7000,-10.5000){\makebox(0,0)[lb]{$1$}}%
\end{picture}}%

\caption{The functions $\psi_{\rho,\rho'}$.}
\end{figure}
Define a probability measure $\mu_{\rho,\rho'}$ on $\DD$ by
$$\mu_{\rho,\rho'}=\frac{1}{\int_\DD\psi_{\rho,\rho'}v_P}\,\,\psi_{\rho,\rho'}\,v_P,$$
where $v_P$ denotes the Poincar\'e volume form. 

Let $(M,\FF)$ be as in Section 1. For any $x\in M$, let $L_x$ be the
leaf through $x$ with the universal cover $\tilde L_x$ identified with
$\DD$. Define a map $j_x:\DD\to M$ as the composite
$$
j_x:\DD\cong \tilde L_x\to L_x\subset M
$$
such that $j_x(0)=x$. Define $\mu_{\rho,\rho'}(x)\in\PPP(M)$ by
$\mu_{\rho,\rho'}(x)=(j_x)_*\mu_{\rho,\rho'}$.
The main result of this section is the following.

\begin{theorem}\label{t3.1}
 If $\mu_{\rho_n,\rho_n'}(x_n)$ converges for some sequences $x_n\in M$,
$\rho_n\to\infty$ and $\rho_n'\to\infty$, then the limit is a harmonic
 measure for $\FF$.
\end{theorem}

To show this, we approximate $\psi_{\rho,\rho'}$ by another
function $\varphi_{\rho,\rho'}$ which is a combination of harmonic
functions.
Let
$A=1/\log\frac{R'}{R}$. Define a function
$\varphi_{\rho,\rho'}:\DD\to[0,\infty)$ by
$$
\varphi_{\rho,\rho'}(z)=\left\{ \begin{array}{ll} 1 & \mbox{\ if \ }\abs{z}\leq R
\\ A\log\frac{R'}{\abs{z}} & \mbox{\ if \ } R\leq\abs{z}\leq R'\\
0 & \mbox{\ if \ } R'\leq\abs{z}<1.\end{array}
\right.$$ 
Define a probability measure $\nu_{\rho,\rho'}$ on $\DD$ by
$$\nu_{\rho,\rho'}=\frac{1}{\int_\DD\varphi_{\rho,\rho'}v_P}\,\,\varphi_{\rho,\rho'}\,v_P,$$
and define $\nu_{\rho,\rho'}(x)=(j_x)_*\nu_{\rho,\rho'}\in\PPP(M)$ just as before.
Theorem \ref{t3.1} reduces to the following two propositions.
Denote by $\Vert\cdot\Vert$ the norm of $\PPP(M)\subset C(M)'$
dual to the sup norm $\Vert\cdot\Vert_\infty$ of the Banach space $C(M)$ of the continuous functions
of $M$.

\begin{proposition}\label{p3.1}
We have $\Vert \mu_{\rho_n,\rho_n'}(x_n)-\nu_{\rho_n,\rho_n'}(x_n)\Vert\to0$
as $\rho_n,\rho_n'\to\infty$.
 \end{proposition}

\begin{proposition}\label{p3.2}
 If $\nu_{\rho_n,\rho_n'}(x_n)$ converges for some sequences $x_n\in M$ and
$\rho_n,\rho_n'\to\infty$, then the limit is a harmonic
 measure for $\FF$.
\end{proposition}

We shall first show Proposition \ref{p3.2}.
For a $C^2$ function $f:M\to\R$, we denote by $\Delta_Pf$
the leafwise Laplacian with respect to the leafwise Poincar\'e metric.
What we have to prove is
that 
$$\int_M\Delta_P f\,\,\nu_{\rho_n,\rho'_n}(x_n)\to 0\ \mbox{ as }\ n\to\infty.$$


Since $j_{x_n}$ is a local isometry onto the leaf, we have
$\Delta_P f\circ j_{x_n}=\Delta_P(f\circ j_{x_n})$.
Rewriting $f\circ j_{x_n}$ as $f$, 
this follows from the following proposition about $\DD$.

\begin{proposition}\label{p}
For any nonzero bounded $C^2$ function $f:\DD\to\R$, we have
$$ 
\frac{\int_\DD\varphi_{\rho,\rho'}\Delta_Pf\,v_P\,\cdot\,{\Vert f\Vert_\infty}^{-1}}
{\int_\DD\varphi_{\rho,\rho'}v_P}\to0
$$
as $\rho,\rho'\to0$ uniformly on $f$
\end{proposition}

\medskip
{\sc Estimate of the numerator:} First notice that
$\Delta_Pf\,v_P=\Delta_Ef\,v_E$ where $E$ stands for Euclidian.
(Both are equal to $dJ^*df$, where $J$ is the almost complex
structure.) \ We need the following Green-Riesz formula.
See \cite{D}, Chapt.I, p.30.

\begin{theorem}
 Let $\Omega$ be a smoothly bounded compact domain in $\R^n$, and let
$\vec n_E$ be the outward unit normal vector at $\partial \Omega$.
Denote by $\sigma_E$ the Euclidian area measure on $\partial\Omega$.
Then for any $C^2$ function $\varphi$ and $f$ defined
on $\R^n$, we have
$$
\int_\Omega(\varphi\Delta_E
 f-f\Delta_E\varphi)v_E=\int_{\partial\Omega}(\varphi\frac{\partial
 f}{\partial\vec{n}_E}- f\frac{\partial\varphi}{\partial\vec n_E})\sigma_E.
$$
\end{theorem}

Let us apply this formula to $\varphi_{\rho,\rho'}$, $f$
 in Proposition \ref{p} and the domains $\DD(R)$, $\DD(R')\setminus
 \DD(R)$.
Remark that
$$\Delta_E\log\frac{R'}{\abs{z}}=0\ \ \mbox{ and }
\frac{\partial}{\partial\vec n_E}(\log\frac{R'}{\abs{z}})=-\frac{1}{\abs{z}}.$$
Computation shows that
$$\int_\DD\varphi_{\rho,\rho'}\Delta_Pf\,v_P
= A(\frac{1}{R'}\int_{\partial\DD(R')}f\sigma_E
- \frac{1}{R}\int_{\partial\DD(R)}f\sigma_E).
$$
This implies that
\begin{equation*}
 \label{ee1}
\abs{\int_\DD\varphi_{\rho,\rho'}\Delta_Pf\,v_P}\,\,\Vert f\Vert_\infty^{-1}\leq
4\pi A\,\,.
\end{equation*}

\medskip
{\sc Estimate of the denominator:} \
We use the following notations.

\begin{notation}
 For $0<R<R'<1$ and positive valued functions $F(R,R')$ and $G(R,R')$,
we write $F\sim G$ if $F/G\to1$ as $R\to1$.
\end{notation}

\begin{lemma}
 \label{l}
(1) If $F_1\sim G_1$ and $F_2\sim G_2$, then $F_1+F_2\sim G_1+G_2$.

(2) If $F\sim G$, then we have
$$
\int_{R}^{R'}F(r,R')dr\sim \int_R^{R'}G(r.R')dr.$$

(3) We have $\log(R'/R)\sim R'-R$.
\end{lemma}

\medskip

Now since
$$v_P=\frac{4dx\wedge dy}{(1-\abs{z}^2)^2}=\frac{4rdr\wedge d\theta}{(1-r^2)^2}$$
in the polar coordinates, we have
$$
\frac{1}{8\pi}\int_\DD\varphi_{\rho,\rho'}v_P
=\int_0^R\frac{rdr}{(1-r^2)^2}+A\int_R^{R'}\log\frac{R'}{r}\,\,
\frac{rdr}{(1-r^2)^2},$$
where
$$
\mbox{the first term\ }
=\frac{R^2}{2(1-R)(1+R)}\sim\frac{1}{4(1-R)}
=\frac{A\log(R'/R)}{4(1-R)}\sim\frac{A(R'-R)}{4(1-R)},$$
and by Lemma \ref{l} (2)
$$
\mbox{the second term\ }\sim A\int_R^{R'}\frac{(R'-r)dr}{4(1-r)^2}
=-\frac{A(R'-R)}{4(1-R)}+\frac{A}{4}\log\frac{1-R}{1-R'}.
$$
Since both terms are positive, we get from Lemma \ref{l} (1),
$$\frac{1}{8\pi}\int_\DD\varphi_{\rho,\rho'}v_P
\sim \frac{A}{4}\log\frac{1-R}{1-R'}\sim \frac{A\rho'}{2},
$$
where the last $\sim$ holds when $\rho'$ is bounded from below
and follows from the formula
$$
\rho'=\frac{1}{2}(\log\frac{1+R'}{1-R'}-\log\frac{1+R}{1-R}).
$$

It follows from the two estimates that
$$
\lim_{R\to1}
\frac{\int_\DD\varphi_{\rho,\rho'}\Delta_Pf\cdot\Vert f\Vert_\infty^{-1}}
{\int_\DD\varphi_{\rho,\rho'}v_P}
\leq \lim_{R\to1}\frac{1}{\rho'}=0,$$
and the convergence is uniform on $f$.
This shows Propositions \ref{p} and \ref{p3.2}. 

\medskip
Finally Proposition \ref{p3.1} follows from the following estimate
$$
 \int_\DD\psi_{\rho,\rho'}v_P\sim\int_\DD\varphi_{\rho,\rho'}v_P,
$$
since $\psi_{\rho,\rho'}\geq\varphi_{\rho,\rho'}$. We have already shown
that 
$$\int_\DD\varphi_{\rho,\rho'}v_P\sim4\pi A\rho'.$$
Analogous (and easier) computation using $A(R'-R)\sim 1$ shows that
$$\int_\DD\psi_{\rho,\rho'}v_P\sim4\pi A\rho'.$$

\section{Further question and example}
It seems that the counterexample in Section 2 is rather special.
There might be more foliations which satisfy the conclusion
of the  theorem of Bonatti and
G\'omez-Mont. To consider this problem, let us recall their proof, which
consists of two steps. In the first step, they
consider general $(M,\FF,g_P))$ as in the beginning of Section 1.
Let $p:\hat M\to M$ be the unit tangent bundle of the foliation $\FF$.
The space $\hat M$ admits a leafwise geodesic flow $\{g^t\}$ and the
leafwise stable horocycle flow $\{h^s\}$, which satisfy
\begin{equation}\label{eq}
g^t\circ h^s\circ g^{-t}= h^{se^{-t}}. 
\end{equation}
Therefore the two flows form a locally free action of the Lie group $B$,
the 2 dimensional nonabelian Lie group. Given a leafwise submersed
$\rho$-disk $B_\rho(z)$ of $\FF$ (See Section 1), they considered the lift
$\sigma:B_\rho(z)\setminus\{z\}\to\hat M$ by the radial unit vector
fields, and showed that the limit
$\displaystyle\lim_{\rho_n\to\infty}\sigma_*\beta_{\rho_n}(z_n)$ is
$h^s$-invariant, if it exists. This part is true for any $(M,\FF,g_P)$.

On the other hand, Yuri Bakhtin and Matilde Mart\'inez \cite{BM} showed that the map
$p_*:\PPP(\hat M)\to\PPP(M)$
between the space of the probability measures gives a bijection from
the subset of the $B$-invariant measures on $\hat M$ to the subset
of the harmonic measures on $M$.

Now assume that the horocycle flow $\{h^s\}$ is uniquely ergodic. Then
the unique invariant measure $\hat\mu$ is also $g^t$-invariant by
(\ref{eq}), and thus $p_*\hat\mu$ is a unique harmonic measure of $(M,\FF,g)$.
In the second step, Bonatti and
Gom\'ez-Mont showed the unique ergodicity of the horocycle flow
$\{h^s\}$ for
foliations in Theorem 1.1. It is plausible to expect that there are
more foliations with this property. We shall raise one example.

\begin{example} \label{ex1}
 Let $G$ be an arbitrary connected unimodular Lie group, and let 
$\Gamma\subset PSL(2,\R)\times G$ be a cocompact lattice such that
$p_2(\Gamma)$ is dense in $G$, where $p_2$ is the projection onto the
 second factor. Then the manifold $M=\Gamma\setminus(\HH\times G)$ admits
a horizontal foliation 
$\FF=\Gamma\setminus\{\HH\times\{g\}\}$ by hyperbolic surfaces. 
\end{example}

The unit tangent bundle $\hat M$ 
of the above foliation $\FF$ is identified
with $\Gamma\setminus(PSL(2,\R)\times G)$.
According to Marina Ratner \cite{R}, any probability
measure $\hat\mu$ invariant by the leafwise stable
horocycle flow is {\em algebraic} in
the following sense. For any $x$ in the support of $\hat\mu$, there is a closed subgroup
$H\subset PSL(2,\R)\times G$ such that the closure of the horocycle orbit of $x$
is $x\cdot H$ and that $\hat\mu=x_*m$, where $m$ is the normalized
Haar measure of
$(g^{-1}\Gamma g\cap H)\setminus H$ and $x=\Gamma g$.  

On the other hand, it is shown by Fernando Alcalde Cuesta and Fran\c{c}oise Dal'bo  \cite{ACD}
that the leafwise stable horocycle flow is minimal. Therefore $\hat\mu$
is the Haar measure of $\Gamma\setminus(PSL(2,\R)\times G)$ and is unique.
In conclusion the foliation in Example 4.1  is equidistributed,
i.e.\ satisfies the conclusion of Theorem 1.1.

\end{document}